\documentclass[a4paper]{amsart}

\usepackage{amsfonts}
\usepackage{verbatim}
\usepackage{amssymb}
\usepackage{amscd}
\usepackage{amsmath}
\usepackage{latexsym}
\usepackage{framed, color}

\usepackage[dvips]{graphicx}
\usepackage{psfrag}

\newtheorem{thm}{Theorem}[section]
\newtheorem{cor}[thm]{Corollary}
\newtheorem{prop}[thm]{Proposition}
\newtheorem{lem}[thm]{Lemma}

\theoremstyle{definition}
\newtheorem{dfn}[thm]{Definition}

\theoremstyle{remark}
\newtheorem{rem}[thm]{Remark}
\newtheorem{problem}[thm]{Problem}

\makeatletter
 
 \@addtoreset{equation}{section}
\makeatother

\newcommand{\R}{\mathbb{R}}
\newcommand{\C}{\mathbb{C}}

\newcommand{\rank}{\mathop{\mathrm{rank}}\nolimits}

\def\spmapright#1{\smash{%
 \mathop{\hbox to 1.3cm{\rightarrowfill}}
  \limits^{#1}}}
\def\spmapleft#1{\smash{%
 \mathop{\hbox to 1.3cm{\leftarrowfill}}
  \limits^{#1}}}

\title[Differentiable maps on links of complex isolated 
hypersurface singularities]
{Differentiable maps on links of complex isolated hypersurface singularities}
\dedicatory{Dedicated to Professor Goo Ishikawa
on the occasion of his retirement}
\author{Osamu Saeki and Shuntaro Sakurai} 
\address{Institute of Mathematics for Industry,
Kyushu University,
Motooka 744, Nishi-ku, Fukuoka 819-0395, Japan}

\address{Daiwa House Industry}

\date{\today}
\keywords{complex isolated hypersurface singularity,
singularity link,
fold map, round fold map, graph manifold}

\subjclass[2010]{Primary 57R45;
Secondary 32S50,
58K30,
58K15}

\begin{document}
\begin{abstract}
We consider links of complex isolated
hypersurface singularities in $\C^{n+1}$
and study differentiable
maps defined by restricting holomorphic
functions to the links. We give an explicit
example in which such a restriction gives a fold map
into the plane $\C = \R^2$
whose singular value set consists
of concentric circles.
\end{abstract}

\maketitle

\section{Introduction}\label{section1}

Let us consider an isolated singularity
of a complex surface in $\C^3$.
Its link is a closed connected orientable
$3$--dimensional manifold \cite{Milnor68}, and it is known to
be a graph manifold \cite{Mumford, N}.
Here, a \emph{graph manifold} is a closed orientable
$3$--dimensional manifold obtained from
some $S^1$--bundles over compact surfaces with boundary
attached along their torus boundaries.

On the other hand, graph manifolds can
be characterized by using their $C^\infty$ stable
maps into $\R^2$ \cite{Sa2, Sa1}.
Recall that smooth generic maps of manifolds of dimension $\geq 2$
into $\R^2$
have folds and cusps as their singularities \cite{L0, T, Wh},
and $C^\infty$ stable maps, in general, have such
singularities.
In \cite{KS2}, it has been proved
that a closed orientable $3$--dimensional manifold
is a graph manifold if and only if it admits a round
fold map into $\R^2$, where a round fold map
is a smooth map that has only folds as its singularities
and its restriction to the singular point set
is an embedding onto a family of concentric circles
in $\R^2$ (for details, see \S\ref{section2}).

Then, a natural question arises: can we construct
an explicit round fold map into $\R^2 = \C$ on
a link of an isolated surface singularity?
Recall that, by a celebrated result due to Mather
\cite{Mather}, if we restrict a real linear projection
$\C^3 \to \R^2$ to the link of a surface singularity,
then generically it is a $C^\infty$ stable map.
However, it has, in general, cusps and its
singular value set might be very complicated.
So, an explicit construction seems to be a challenging
problem.

In this paper, we consider the link of the
hypersurface singularity in $\C^{n+1}$
defined by the Brieskorn--Pham type polynomial
$$f(z_1, z_2, \ldots, z_{n+1}) =
z_1^2 + z_2^2 + \cdots + z_{n+1}^2$$
and construct an explicit round fold
map on its link by restricting a complex
linear function on $\C^{n+1}$. Although the example
is very simple and the proof consists of tedious,
but elementary computations, this seems to be
a first example in which a simple generic map
is constructed explicitly.

The paper is organized as follows. In \S\ref{section2},
we first formulate our motivating problem, then
we give some definitions and preliminary results
for determining the singular point set
of the relevant smooth map into the plane.
In \S\ref{section3}, we consider the explicit
example as mentioned above and show that
the restriction of a complex linear function to
the hypersurface singularity link
gives a round fold map into the plane and
we also identify the indices of the fold loci.

Throughout the paper, 
all manifolds and maps between them are smooth
of class $C^\infty$ unless otherwise specified. 
For a map $\varphi : M \to N$ between manifolds,
we denote by $S(\varphi)$ the set of points in $M$
where the differential of $f$ does not have
maximal rank $\min \{\dim{M}, \dim{N}\}$.
The symbol ``$\cong$'' denotes a diffeomorphism between
smooth manifolds.

\section{Preliminaries}\label{section2}

\subsection{Problem}\label{subsec2-1}

\begin{figure}[t]
\centering
\psfrag{C}{$\C^{n+1}$}
\psfrag{S}{$S^{2n+1}_\varepsilon$}
\psfrag{K}{$K_f$}
\psfrag{0}{$\mathbf{0}$}
\psfrag{f}{$f^{-1}(0)$}
\includegraphics[width=0.95\linewidth,height=0.3\textheight,
keepaspectratio]{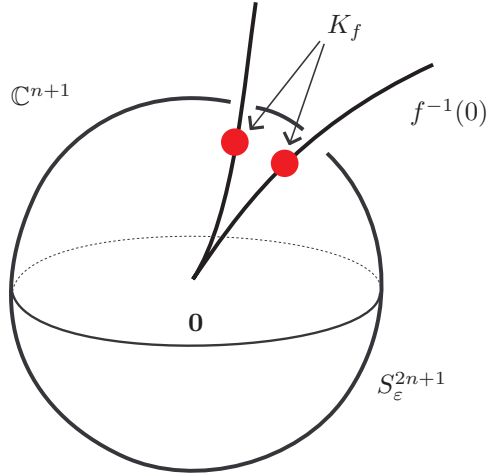}
\caption{Link of an isolated hypersurface singularity}
\label{fig1}
\end{figure}

Let $f : (\C^{n+1}, \mathbf{0}) \to (\C, 0)$, $n \geq 1$, be
a holomorphic function germ with an isolated
critical point at the origin.
In the following,
we abuse the notation and denote a representative
of the germ by the same letter.
It is known that for a sufficiently small $\varepsilon > 0$,
$$K_f = f^{-1}(0) \cap S^{2n+1}_\varepsilon$$
is a smooth closed orientable $(n-2)$--connected
$(2n-1)$--dimensional manifold, called
the \emph{link} of $f$, where
$S^{2n+1}_\varepsilon$ is the sphere of radius $\varepsilon$
centered at the origin in $\C^{n+1}$ \cite{Milnor68}
(see Fig.~\ref{fig1}).
Note also that the diffeomorphism type of $K_f$
and also its isotopy class in $S^{2n+1}_\varepsilon$
do not depend on the choice of $\varepsilon$
as long as it is sufficiently small.

Now let $g : U \to \C$ be an arbitrary
holomorphic function defined on an open neighborhood
of the origin $\mathbf{0}$ in $\C^{n+1}$, and set 
$h = g|_{K_f} : K_f \to \C = \R^2$, where
we take $\varepsilon > 0$ so small that $K_f$
is contained in $U$.
The problem that we consider in the
present paper is the following.

\begin{problem}\label{prob1}
Study the smooth map $h$ from
a global singularity theoretical viewpoint.
\end{problem}

\subsection{Round fold maps}\label{subsec2-2}

We will see that in a special case, such a map
$h$ as in Problem~\ref{prob1} is a round fold map 
as defined below.

\begin{dfn}
Let $M$ be a smooth $m$--dimensional manifold and
let $\varphi : M \to \R^p$, $m \geq p \geq 1$, be a smooth map.
A singular point $q \in M$ of $\varphi$ is called
a \emph{fold} if $\varphi$ can be written
as in the following normal form with respect to some local
coordinates around $q$ and $\varphi(q)$:
$$(x_1, x_2, \ldots, x_m) \mapsto (x_1, x_2, \ldots, x_{p-1}, \pm x_p^2 \pm
x_{p+1}^2 \pm \cdots \pm x_m^2).$$
Let $\lambda$ denote the number of negative signs
appearing in the last quadratic form as above.
Then, the integer $\min\{\lambda, m-p+1-\lambda\}$ is well-defined
and is called the \emph{absolute index} of the fold $q$.
A fold is \emph{definite} if its absolute index is equal to zero,
otherwise it is \emph{indefinite}.
If $\varphi$ has only fold as its singularities, then it
is called a \emph{fold map} \cite{E}.
\end{dfn}

Note that if $\varphi$ is a fold map, then
its singular point set $S(\varphi)$ is a
$(p-1)$--dimensional closed submanifold of $M$
and the restriction $\varphi|_{S(\varphi)} :
S(\varphi) \to \R^p$
is a codimension one immersion.

\begin{dfn}
A smooth map $\varphi : M \to \R^p$
is called a \emph{round fold map}
if it is a fold map, each connected component
of $S(\varphi)$ is diffeomorphic to the
standard $(p-1)$--dimensional sphere, and
the map $\varphi|_{S(\varphi)}$
is an embedding whose image is
isotopic to a family of concentric spheres.
For details, see \cite{Ki0, Ki1, Ki2}.
\end{dfn}

\subsection{Inner products}\label{subsec2-3}

In the following, for 
$$\mathbf{u} = (u_1, u_2, \ldots, u_{n+1}), \,
\mathbf{v} = (v_1, v_2, \ldots, v_{n+1}) \in \C^{n+1},$$
their \emph{Hermitian inner product} is defined as
$$\langle \mathbf{u}, \mathbf{v} \rangle
= \sum_{j=1}^{n+1} u_j \bar{v}_j \in \C.$$
On the other hand, if we regard
$\mathbf{u}, \mathbf{v} \in \C^{n+1} = \R^{2n+2}$
and set $u_j = u_j^0 + \sqrt{-1} u_j^1$, $v_j = v_j^0 + \sqrt{-1} v_j^1$
with $u_j^0, u_j^1, v_j^0, v_j^1 \in \R$,
then their \emph{real inner product} is defined as
$$(\mathbf{u}, \mathbf{v})_\R = 
\sum_{j=1}^{n+1} (u_j^0 v_j^0 +
u_j^1 v_j^1) \in \R.$$
Then, the following is easily verified.

\begin{lem}\label{lemma1.1} We always have
$\displaystyle{(\mathbf{u}, \mathbf{v})_\R = \mathrm{Re} \langle
\mathbf{u}, \mathbf{v} \rangle.}$
\end{lem}

Let $f : U \to \C$ be a holomorphic function defined
on an open set $U \subset \C^{n+1}$.
For $\mathbf{z}_0 \in U$, we set
$$\overline{\mathbf{grad}}\,f(\mathbf{z}_0) := \left(
\overline{\frac{\partial f}{\partial z_1}(\mathbf{z}_0)},
\overline{\frac{\partial f}{\partial z_2}(\mathbf{z}_0)},
\ldots, 
\overline{\frac{\partial f}{\partial z_{n+1}}(\mathbf{z}_0)}
\right).$$
We warn the reader that this notation is different
from that used in \cite{Milnor68}.

\begin{lem}\label{lem2}
Let $y \in \C$ be a regular value of $f$
and set $M = f^{-1}(y) \subset \C^{n+1}$.
Then, the complex tangent space $T_{\mathbf{z}_0}M$ of 
the complex manifold $M$ at
$\mathbf{z}_0 \in M$ can be identified with the
$\C$--vector space
$$\overline{\mathbf{grad}}\,f(\mathbf{z}_0)^{\perp_\C}
= \{\mathbf{v} \in \C^{n+1}\,|\, \langle
\mathbf{v}, \overline{\mathbf{grad}}\,f(\mathbf{z}_0)
\rangle = 0\}.$$
\end{lem}

\begin{proof}
Let $\mathbf{p}(t) = (p_1(t), p_2(t), \ldots, p_{n+1}(t))$ 
$(t \in (-\delta, \delta), \delta > 0)$
be a $C^\infty$ curve in $M$ with $\mathbf{p}(0) =
\mathbf{z}_0$. Since
$f(\mathbf{p}(t)) = y$ is constant, we have
\begin{eqnarray*}
\left.\frac{d}{dt}f(\mathbf{p}(t))\right|_{t=0}
& = & \sum_{j=1}^{n+1} \frac{\partial f}{\partial z_j}(\mathbf{z}_0)
\frac{d p_j}{dt}(0) \\
& = & \langle \mathbf{p}'(0), \overline{\mathbf{grad}}\,
f (\mathbf{z}_0) \rangle \, = \, 0.
\end{eqnarray*}
Therefore, we have
$$T_{\mathbf{z_0}}M \subset
\overline{\mathbf{grad}}\,f(\mathbf{z_0})^{\perp_\C}.$$
Since the two vector spaces over $\C$ have the same dimension, we get the
desired result.
\end{proof}

\subsection{Singular point set}\label{subsec2-4}

Now, let us go back to the situation of Subsection~\ref{subsec2-1}:
$f : (\C^{n+1}, \mathbf{0}) \to (\C, 0)$
is a holomorphic function germ with an isolated
critical point at the origin, 
$K = K_f = f^{-1}(0) \cap S^{2n+1}_\varepsilon$
is the link associated with $f$ for a sufficiently small
$\varepsilon > 0$, $g : U \to \C$
is a holomorphic function defined on an open
neighborhood $U$ of the origin $\mathbf{0}$ in $C^{n+1}$, 
and consider the smooth map $h : K \to \C = \R^2$ defined
by $h = g|_K$.
Note that $h$ is a smooth map
of a smooth closed $(2n-1)$--dimensional
manifold $K$ into $\R^2$. In the following, we assume $n \geq 2$.
We denote by $S(h)$ $(\subset K)$ the
set of singular points of $h$.

\begin{lem}\label{lemma2.3}
For $\mathbf{z}_0 \in K$, we have
$\mathbf{z}_0 \in S(h)$ if and only if
the three vectors
$\overline{\mathbf{grad}}\, f(\mathbf{z}_0)$,
$\overline{\mathbf{grad}}\, g(\mathbf{z}_0)$, and
$\mathbf{z}_0$ in $\C^{n+1}$ are linearly dependent over $\C$.
\end{lem}

\begin{proof}
For $\mathbf{z}_0 \in K$, 
since $f^{-1}(0)$ and $S^{2n+1}_\varepsilon$ 
intersect transversely \cite{Milnor68}, by Lemma~\ref{lem2}
and its proof,
we have
$$T_{\mathbf{z}_0}K =
\overline{\mathbf{grad}}\, f (\mathbf{z}_0)^{\perp_\C}
\cap \mathbf{z}_0^{\perp_\R},$$
where
$$T_{\mathbf{z}_0}S^{2n+1}_\varepsilon =
\mathbf{z}_0^{\perp_\R} =
\{\mathbf{v} \in \C^{n+1}\,|\, (\mathbf{v}, \mathbf{z}_0)_\R = 0\}.$$
Then, we see easily the following equivalences:
\begin{eqnarray*}
& & \mathbf{z}_0 \in S(h) = S(g|_K) \\
& \Longleftrightarrow &
\dim_{\R}{\mathrm{Ker}\,(dg_{\mathbf{z}_0}|_{T_{\mathbf{z}_0K}})}
\geq 2n-2 \\
& \Longleftrightarrow &
\dim_{\R} ((\overline{\mathbf{grad}}\, f (\mathbf{z}_0)^{\perp_\C}
\cap \mathbf{z_0}^{\perp_\R})
\cap \overline{\mathbf{grad}}\, g(\mathbf{z}_0)^{\perp_\C}) \geq 2n-2 \\
& \Longleftrightarrow &
\dim_{\R} ((\overline{\mathbf{grad}}\, f (\mathbf{z}_0)^{\perp_\C}
\cap \overline{\mathbf{grad}}\, g(\mathbf{z}_0)^{\perp_\C}) 
\cap \mathbf{z_0}^{\perp_\R})
\geq 2n-2.
\end{eqnarray*}
Note that $\overline{\mathbf{grad}}\, f (\mathbf{z}_0)$
is non-zero by our assumption, while 
$\overline{\mathbf{grad}}\, g(\mathbf{z}_0)$
might be zero. Then, the above condition,
in turn, is equivalent to (\ref{eq1}) or (\ref{eq2})
as described below:
\begin{eqnarray}
& & \mbox{$\overline{\mathbf{grad}}\, f(\mathbf{z}_0)$ and
$\overline{\mathbf{grad}}\, g(\mathbf{z}_0)$ are linearly
dependent over $\C$, or} 
\label{eq1}
\\
& & 
V := \overline{\mathbf{grad}}\, f (\mathbf{z}_0)^{\perp_\C}
\cap \overline{\mathbf{grad}}\, g(\mathbf{z}_0)^{\perp_\C} 
\subset \mathbf{z}_0^{\perp_\R}.
\label{eq2}
\end{eqnarray}
If (\ref{eq2}) holds, then $\forall \mathbf{v} \in V$, we have
$0 = (\mathbf{v}, \mathbf{z}_0)_\R = \mathrm{Re}\langle
\mathbf{v}, \mathbf{z}_0 \rangle$.
As $V$ is a vector space over $\C$ and we have
$0 = (-\sqrt{-1} \mathbf{v}, \mathbf{z}_0)_{\R} = 
\mathrm{Re}\langle
-\sqrt{-1} \mathbf{v}, \mathbf{z}_0 \rangle = 
\mathrm{Im}\langle
\mathbf{v}, \mathbf{z}_0 \rangle$.
Hence, we have $\langle
\mathbf{v}, \mathbf{z}_0 \rangle = 0$ and hence
$\mathbf{v} \in \mathbf{z}_0^{\perp_\C}$.
Therefore,
(\ref{eq2}) implies $V \subset \mathbf{z}_0^{\perp_\C}$.
Note that the converse obviously holds.
Furthermore, $V \subset \mathbf{z}_0^{\perp_\C}$ is equivalent to that
$\mathbf{z}_0$ is a linear combination of
$\overline{\mathbf{grad}}\, f(\mathbf{z}_0)$ and 
$\overline{\mathbf{grad}}\, g(\mathbf{z}_0)$ over $\C$.
This completes the proof.
\end{proof}

\begin{rem}
Kamiya \cite[Theorem~1]{Kamiya} obtains a general result in the
real $C^\infty$ setting similar to the above lemma.
The above lemma can be considered to be a refined
version of Kamiya's result in our situation
involving links of complex isolated singularities.
\end{rem}

When $n=2$, the condition described in Lemma~\ref{lemma2.3}
is equivalent to that the $3 \times 3$ complex matrix
consisting of the three column vectors
${}^T\overline{\mathbf{grad}}\, f(\mathbf{z}_0)$,
${}^T\overline{\mathbf{grad}}\, g(\mathbf{z}_0)$, and
${}^T\mathbf{z}_0$ satisfies
$$\det 
\begin{pmatrix} {}^T\overline{\mathbf{grad}}\, f(\mathbf{z}_0) &
{}^T\overline{\mathbf{grad}}\, g(\mathbf{z}_0) &
{}^T\mathbf{z}_0
\end{pmatrix}
= 0,$$
where ``$T$'' means the transpose.
Therefore, we have the following.

\begin{prop} 
For $n=2$, the singular point set $S(h)$
of the smooth map $g|_K = h : K \to \R^2$
is given by
$$S(h) = \{\mathbf{z} \in K\,|\, 
\det 
\begin{pmatrix} {}^T\overline{\mathbf{grad}}\, f(\mathbf{z}) &
{}^T\overline{\mathbf{grad}}\, g(\mathbf{z}) &
{}^T\mathbf{z}
\end{pmatrix}
= 0.\}$$
\end{prop}

\section{Explicit construction}\label{section3}

Let us consider the Brieskorn--Pham type polynomial
$$f(\mathbf{z}) = z_1^2 + z_2^2 + \cdots + z_{n+1}^2,$$
where $\mathbf{z} = (z_1, z_2, \ldots, z_{n+1})$.
For the moment, we assume $n \geq 2$.
Setting
$\varepsilon = 1$, its link is given by
$$K_f = \{\mathbf{z} \in \C^{n+1}\,|\, f(\mathbf{z}) = 0, ||\mathbf{z}||^2 = 1\}.$$
Furthermore, let us consider the holomorphic function
$g : \C^{n+1} \to \C$ defined by 
$$g(\mathbf{z}) = z_1 + \frac{\sqrt{-1}}{2} z_2.$$
In fact, the choice of $g$ was made so that we get
the desired result as follows; however, according to the
authors' intuition, to find
such an explicit example is not an easy task.

\begin{rem}
It is known that the $(2n-1)$--dimensional manifold
$K_f$ is diffeomorphic to the total space
of the unit tangent sphere bundle over $S^n$.
When $n=2$, it is the total space of the $S^1$--bundle
over $S^2$ with Euler number $2$, and is diffeomorphic to
the real projective space $\R P^3$, or the lens space $L(2, 1)$.
For general $n$, the singularity of $f$ 
at the origin is called the simple singularity of
type $A_1$, and in a certain sense, it
can be considered to be the simplest complex singularity.
\end{rem}

The main result of this paper is the following.

\begin{thm}\label{thm1}
The smooth map
$h = g|_{K_f} : K_f \to \C = \R^2$ is a
round fold map. The singular point set
$S(h)$ consists of two circles; one consists of
definite folds and the other consists of
indefinite folds of absolute index $n-1$.
\end{thm}

\begin{figure}[t]
\centering
\psfrag{C}{$\C \supset h(S(h))$}
\psfrag{r}{$\sqrt{2}/4$}
\psfrag{s}{$3\sqrt{2}/4$}
\psfrag{R}{$\mathrm{Re}$}
\psfrag{0}{$\mathbf{0}$}
\psfrag{I}{$\mathrm{Im}$}
\includegraphics[width=0.95\linewidth,height=0.5\textheight,
keepaspectratio]{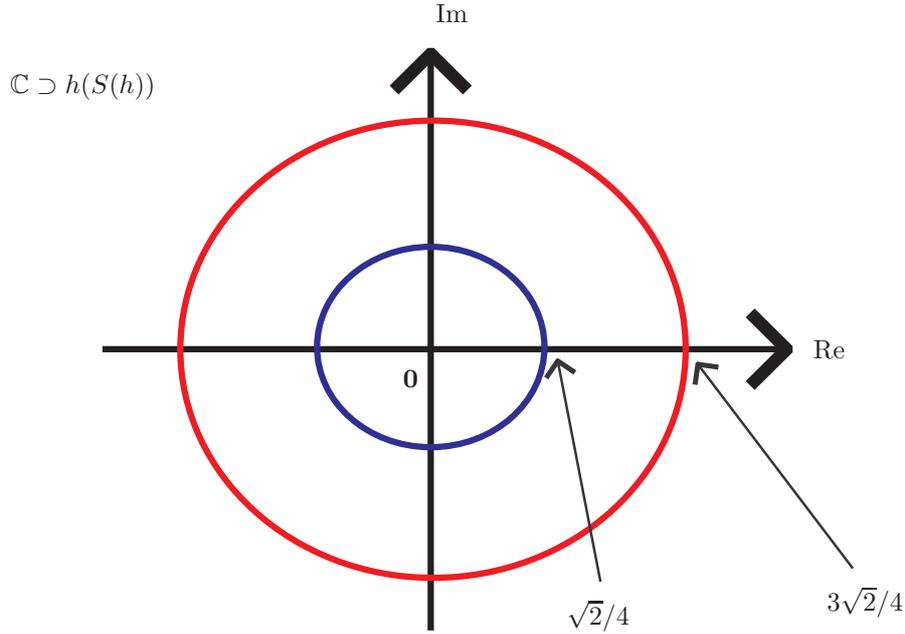}
\caption{Singular value set $h(S(h))$ of $h$}
\label{fig2}
\end{figure}

For the image $h(S(h))$ of the singular point set, see Fig.~\ref{fig2}.

\begin{proof}
First, by Lemma~\ref{lemma2.3}, for $\mathbf{z} \in K_f$,
we have $\mathbf{z} \in S(h)$
if and only if
\begin{equation}
\rank
\begin{pmatrix} {}^T\overline{\mathbf{grad}}\, f(\mathbf{z}) &
{}^T\overline{\mathbf{grad}}\, g(\mathbf{z}) &
{}^T\mathbf{z}
\end{pmatrix}
\leq 2,
\label{eq:rank}
\end{equation}
where
\begin{eqnarray*}
\overline{\mathbf{grad}}\, f(\mathbf{z}) & = & 2(\bar{z}_1, \bar{z}_2, \ldots,
\bar{z}_{n+1}), \\
\overline{\mathbf{grad}}\, g(\mathbf{z}) & = & (1, -\sqrt{-1}/2, 0, 0, \ldots, 0), \\
\mathbf{z} & = & (z_1, z_2, \ldots, z_{n+1}).
\end{eqnarray*}
Suppose $\mathbf{z} \in S(h)$, then we have, for $j \geq 3$,
$$\det
\begin{pmatrix} 2 \bar{z}_1 & 1 & z_1 \\
2 \bar{z}_2 & -\sqrt{-1}/2 & z_2 \\
2 \bar{z}_j & 0 & z_j
\end{pmatrix}
= 4 \sqrt{-1} \, \mathrm{Im}(z_2 \bar{z}_j) - 2 \, \mathrm{Im}(z_1 \bar{z}_j) = 0,
$$
which implies that $z_1 \bar{z}_j$ and $z_2 \bar{z}_j$ 
are real numbers. Furthermore, for $j \geq 4$, we have
$$\det
\begin{pmatrix} 2 \bar{z}_1 & 1 & z_1 \\
2 \bar{z}_3 & 0 & z_3 \\
2 \bar{z}_j & 0 & z_j
\end{pmatrix}
= 4 \sqrt{-1} \, \mathrm{Im}(z_3 \bar{z}_j) = 0,
$$
which implies that $z_3 \bar{z}_j$ is a real number.
Now suppose that $z_3 \neq 0$. Then, the above conditions imply that
$z_1, z_2, z_4, z_5, \ldots, z_{n+1}$ are all real multiples of $z_3$.
Since $f(\mathbf{z}) = z_1^2 + z_2^2 + \cdots + z_{n+1}^2 = 0$,
this implies that $z_3 = 0$, which is a contradiction. Hence, we must have $z_3 = 0$.
By similar arguments, we see that $z_j = 0$ for all $j \geq 3$.

Consequently, since $\mathbf{z} \in K_f$, we have
$$
S(h) \subset \{(z_1, z_2, 0, \ldots, 0) \in \C^{n+1}\,|\,
z_1 = \pm \sqrt{-1} z_2, |z_1| = |z_2| = 1/\sqrt{2}\}.
$$
On the other hand, if $\mathbf{z}$ belongs to
the set on the right hand side above, then we see easily that (\ref{eq:rank}) holds.
Therefore, we have
$$
S(h) = \{(z_1, z_2, 0, \ldots, 0) \in \C^{n+1}\,|\,
z_1 = \pm \sqrt{-1} z_2, |z_1| = |z_2| = 1/\sqrt{2}\}
\cong S^1 \cup S^1.
$$

Furthermore, as
$$h(z_1, z_2, 0, \ldots, 0) = z_1 + \frac{\sqrt{-1}}{2}z_2 
= \begin{cases} 3 \sqrt{-1} z_2/2, & \\
-\sqrt{-1}z_2/2, & 
\end{cases}
$$
we see that $h|_{S(h)} : S(h) \to \C$
is an embedding and that its image
is a family of two concentric circles as
depicted in Fig.~\ref{fig2}.

Let us now show that $h$ is a fold map.
Set
$L := [0, \infty) \times \{0\} \subset \R^2 = \C$.
Then, 
$Q = h^{-1}(L)$ is a compact $(2n-2)$--dimensional
manifold with non-empty boundary. Let us
consider the smooth function
$\psi = h|_Q : Q \to L = [0, \infty)$.

\begin{lem}\label{lem4}
The function $\psi$ is a Morse function with exactly two
critical points with indices $2n-2$ and $n-1$.
\end{lem}

\begin{proof}
Note that by the above computation about $S(h)$,
we see that the critical point set of $\psi$ 
coincides with
$S(h) \cap Q = \{q, q'\}$, where
\begin{eqnarray*}
q & = & \left(\sqrt{2}/2, -\sqrt{-1}\sqrt{2}/2, 0, 0, \ldots, 0\right), \\
q' & = & \left(\sqrt{2}/2, \sqrt{-1}\sqrt{2}/2, 0, 0, \ldots, 0\right).
\end{eqnarray*}
Let us first calculate the Hessian at the critical point $q$.

We easily see that
$$\mathbf{z}' = (z_3, z_4, \ldots, z_{n+1})$$
can be chosen as local coordinates around $q$ for $Q$.
More precisely,
we consider the parametrization
$$\Phi(\mathbf{z}') = 
\left(\frac{\sqrt{2}}{2}+a(\mathbf{z}'), -\frac{\sqrt{2}}{2}\sqrt{-1} +
b(\textbf{z}')\sqrt{-1}, z_3, \ldots, z_{n+1}\right),$$
where $a(\mathbf{z}')$ and $b(\mathbf{z}')$ are complex
functions of class $C^\infty$ with $a(\mathbf{0}) = b(\mathbf{0}) = 0$
determined by the following:
\begin{eqnarray}
& & \left(\frac{\sqrt{2}}{2}+a(\mathbf{z}')\right)^2 -
\left(-\frac{\sqrt{2}}{2} +
b(\textbf{z}')\right)^2 + z_3^2 + \cdots + z_{n+1}^2 = 0, 
\label{eq:para1}\\
& & \left(\frac{\sqrt{2}}{2}+a(\mathbf{z}')\right)
\left(\frac{\sqrt{2}}{2}+\overline{a(\mathbf{z}')}\right)
+ \left(-\frac{\sqrt{2}}{2} + b(\textbf{z}')\right)
\left(-\frac{\sqrt{2}}{2} + \overline{b(\textbf{z}')}\right) 
\label{eq:para2} \\
& & \qquad + z_3\bar{z}_3 + \cdots + z_{n+1}\bar{z}_{n+1} = 1, 
\nonumber \\
& & \mathrm{Im}\,h(\Phi(\mathbf{z}')) = \mathrm{Im}\,a(\mathbf{z}')
- \frac{1}{2}\mathrm{Im}\,b(\mathbf{z}') = 0.
\label{eq:para3}
\end{eqnarray}

Under this parametrization, by differentiating each of the equations
(\ref{eq:para1}), (\ref{eq:para2}) and (\ref{eq:para3}) twice,
after a tedious but elementary
calculation, we see that, for $3 \leq j \leq n+1$,
$$
\begin{pmatrix}
\frac{\partial^2 h}{\partial z_j^2}(q) & 
\frac{\partial^2 h}{\partial z_j \partial \bar{z}_j}(q) \\
\frac{\partial^2 h}{\partial \bar{z}_j \partial z_j}(q) &
\frac{\partial^2 h}{\partial \bar{z}_j^2}(q) 
\end{pmatrix}
= 
\begin{pmatrix}
- 1/4 & 
- 3/4 \\
- 3/4 &
- 1/4
\end{pmatrix}.
$$
Furthermore, all the other entries of the Hessian of $h$ with respect to
$$(z_3, \bar{z}_3, z_4, \bar{z}_4, \ldots, z_{n+1}, \bar{z}_{n+1})$$
are zero.

Then, for 
$$x_j = \frac{z_j + \bar{z}_j}{2},\, 
y_j = \frac{z_j - \bar{z}_j}{2 \sqrt{-1}},$$
we have
$$
\begin{pmatrix}
\frac{\partial^2 h}{\partial x_j^2}(q) & 
\frac{\partial^2 h}{\partial x_j \partial y_j}(q) \\
\frac{\partial^2 h}{\partial y_j \partial x_j}(q) &
\frac{\partial^2 h}{\partial y_j^2}(q) 
\end{pmatrix}
= 
\begin{pmatrix}
- 2 & 0 \\
0 & -1
\end{pmatrix}.
$$
Therefore, the
Hessian of $h$ with respect to
$(x_3, y_3, x_4, y_4, \ldots, x_{n+1}, y_{n+1})$
is nondegenerate and it has $2n-2$ negative
eigenvalues.

For the critical point
$q' = \left(\sqrt{2}/2, \sqrt{-1}\sqrt{2}/2, 0, 0, \ldots, 0\right)
\in S(h) \cap Q$, we can similarly calculate the
Hessian and in this case, the number of negative eigenvalues of
each $(2 \times 2)$--matrix is equal to $1$ and hence
the index of the critical point is equal to $n-1$.
This completes the proof of Lemma~\ref{lem4}.
\end{proof}

Now, we see that the Lie group
$S^1 = \{\alpha \in \C\,|\, |\alpha| = 1\}$ acts
differentiably
on $K_f$ and on $\C$ from left by the
natural scalar multiplication
(recall that 
$f$ is a homogeneous polynomial):
\begin{eqnarray*}
& & \alpha \in S^1, (z_1, z_2, \ldots, z_{n+1}) \in K_f, z \in \C, \\
& & \alpha \cdot (z_1, z_2, \ldots, z_{n+1}) = (\alpha z_1,
\alpha z_2, \ldots, \alpha z_{n+1}), \\
& & 
\alpha \cdot z = \alpha z.
\end{eqnarray*}
For the action on $\C$, refer to Fig.~\ref{fig3}.

\begin{figure}[t]
\centering
\psfrag{L}{$L$}
\psfrag{S}{$S^1$--action}
\psfrag{R}{$\mathrm{Re}$}
\psfrag{0}{$\mathbf{0}$}
\psfrag{I}{$\mathrm{Im}$}
\includegraphics[width=0.95\linewidth,height=0.5\textheight,
keepaspectratio]{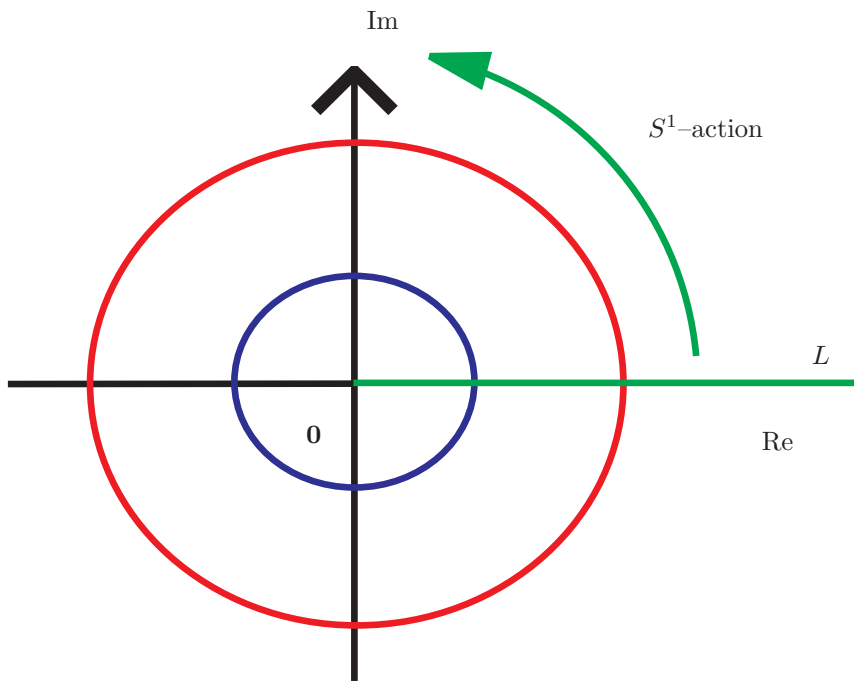}
\caption{$S^1$--action}
\label{fig3}
\end{figure}

As $h$ is a restriction of a linear function,
$h$ is equivariant with respect to the $S^1$--actions:
$$h(\alpha \cdot \mathbf{z}) = \alpha \cdot h(\mathbf{z}).$$
In particular, the singular point set $S(h)$ is $S^1$--invariant.

Furthermore, by the $S^1$--action, $L$ sweeps out $\C$ and
it covers all the singular values.
This, together with the fact that $h|_{S(h)}$
is an embedding, implies that $h$ is a fold map.

Since its singular value set consists of two concentric
circles, it is a round fold map.
This completes the proof of Theorem~\ref{thm1}.
\end{proof}

\begin{rem}
When $n=1$,
$K_f \subset S^3$ is, in fact, a Hopf link.
In this case, 
$h = g|_{K_f}$ is an embedding into $\C$ and
its image is as depicted in Fig.~\ref{fig2}.
\end{rem}

As a corollary of our main theorem, we have
the following, immediately.

\begin{cor}
Let $\eta : \C \to \R$ be an arbitrary nonzero
real linear function.
Then,
$\eta \circ h : K_f \to \R$
is a Morse function with exactly four critical
points of indices
$0$, $n-1$, $n$, and $2n-1$.
\end{cor}

The above corollary can be proved by using
some results concerning the critical
points of such a composite function, obtained in \cite{Fukuda}.

\section*{Acknowledgment}\label{ack}
The authors would like to express their sincere gratitude
to Professor Goo Ishikawa for his constant encouragement.
Some of the results in this paper are contained
in the master thesis of the second author \cite{Sakurai}: the authors would
like to thank the members of the regular seminar, especially
Dr.~Naoki Kitazawa, for their useful discussions.
This work has been supported in part by JSPS KAKENHI Grant Numbers 
JP22K18267, JP23H05437. 


\end{document}